\documentclass[draft]{amsart}
\usepackage{amssymb}

\newcommand\dirlim{\mathop{\varinjlim}\limits}
\newcommand\Zset{\mathbb {Z}}
\newcommand\Qset{\mathbb {Q}}
\newcommand\Cset{\mathbb {C}}
\newcommand\Rset{\mathbb {R}}
\newcommand\im{\mathrm{im}}
\newcommand\tr{\mathrm{tr}}
\newcommand\lge{l^2(G)}

\newcommand\ef{{\mathcal F}}
\newcommand\te{{\mathcal T}}
\newcommand\ce{{\mathcal C}}
\newcommand\vng{{\mathcal N}G}
\newcommand\vnh{{\mathcal N}H}
\newcommand\ug{{\mathcal U}G}

\newcommand\NZ{{\mathcal N}\Zset}
\newcommand\UZ{{\mathcal U}\Zset}
\newcommand\A{{\mathcal A}}
\newcommand\U{{\mathcal U}}
\newcommand\B{{\mathcal B}}
\newcommand\la{l^2(\A)}
\newcommand\lb{l^2(\B)}
\newcommand\T{\mathrm{{\bf T}}}
\newcommand\smallt{\mathrm{{\bf t}}}
\newcommand\bigP{\mathrm{{\bf P}}}
\newcommand\p{\mathrm{{\bf p}}}
\newcommand\unb{\mathrm{{\bf u}}}
\newcommand\bnd{\mathrm{{\bf b}}}
\newcommand\cl{\mathrm{cl}}
\newcommand\tor{\mathrm{Tor}}
\newcommand\homo{\mathrm{Hom}}

\newtheorem{thm}{Theorem}[section]
\newtheorem{lem}{Lemma}[section]
\newtheorem{cor}{Corollary}[section]
\newtheorem{prop}{Proposition}[section]
\newtheorem{defn}{Definition}[section]
\newtheorem{exmp}{Example}[section]

\begin{document}

\title{TORSION THEORIES FOR FINITE VON NEUMANN ALGEBRAS}

\author{Lia Va\v s}

\address{Department of Mathematics, Physics and Computer Science,
University of the Sciences in Philadelphia, 600 S. 43rd St.,
Philadelphia, PA 19104}

\email{l.vas@usip.edu}

\thanks{Part of the results are obtained during the time the author
was at the University of Maryland, College Park. The author was
supported by NSF grant DMS9971648 at that time.}



\begin{abstract}
The study of modules over a finite von Neumann algebra $\A$ can be
advanced by the use of torsion theories. In this work, some
torsion theories for $\A$ are presented, compared and studied. In
particular, we prove that the torsion theory $(\T,\bigP)$ (in
which a module is torsion if it is zero-dimensional) is equal to
both Lambek and Goldie torsion theories for $\A$.

Using torsion theories, we describe the injective envelope of a
finitely generated projective $\A$-module and the inverse of the
isomorphism $K_0(\A)\rightarrow K_0(\U),$ where $\U$ is the
algebra of affiliated operators of $\A.$ Then, the formula for
computing the capacity of a finitely generated module is obtained.
Lastly, we study the behavior of the torsion and torsion-free
classes when passing from a subalgebra $\B$ of a finite von
Neumann algebra $\A$ to $\A$. With these results, we prove that
the capacity is invariant under the induction of a $\B$-module.
\end{abstract}

\maketitle

\section{INTRODUCTION}

Recently there has been an increased interest in the subject of
group von Neumann algebras. For a thorough survey of group von
Neumann algebras, their applications to various fields of
mathematics and list of open problems the interested reader should
check \cite{Lu5}.

One of the reasons for this growing interest is that a group von
Neumann algebra $\vng$ comes equipped with a faithful and normal
trace that enables us to define the dimension of an arbitrary
$\vng$-module. The dimension allows us to consider topological
invariants of a $G$-space in cases when ordinary invariants cannot
be calculated.

Moreover, $\vng$ mimics the ring $\Zset$ in such a way that every
finitely generated module is a direct sum of a torsion and
torsion-free part. The dimension faithfully measures the
torsion-free part while another $L^2$-invariant, the
Novikov-Shubin's, measures the torsion part. Although not without
zero-divisors, a group von Neumann algebra is semihereditary
(i.e., every finitely generated submodule of a projective module
is projective) and an Ore ring. The fact that $\vng$ is an Ore
ring allows us to define the classical ring of quotients, denoted
$\ug$. Besides this algebraic definition, it turns out that,
within the operator theory, $\ug$ can be defined as the algebra of
affiliated operators.

As a ring, $\ug$ keeps all the properties that $\vng$ has and
posses some more. In the analogy that $\vng$ is like $\Zset,$
$\ug$ plays the role of $\Qset.$ Thus, $\ug$ is also a good
candidate for a coefficients ring when working with a $G$-space if
one is not interested in the information that gets lost by the
transfer from $\vng$ to $\ug$ (faithfully measured by the
Novikov-Shubin invariant).

Many results on group von Neumann algebras are obtained by
studying a more general class of von Neumann algebras, the class
of finite von Neumann algebras. Every finite von Neumann algebra
has a normal and faithful trace and, as a consequence, has all of
the properties mentioned above for a group von Neumann algebra. In
Section \ref{FiniteVNA}, we define a finite von Neumann algebra
$\A$, the dimension of $\A$-module and the algebra of affiliated
operators of $\A,$ and list some results on these notions that we
shall use further on.

A finite von Neumann algebra $\A$ has the property that every
finitely generated module is a direct sum of a torsion and a
torsion-free module. However, it turns out that there exists more
than just one suitable candidate when it comes to defining torsion
and torsion-free modules. To clarify the situation, we demonstrate
that the notion of a torsion theory of a ring is a good framework
for the better understanding of the structure of $\A$-modules. In
Section \ref{TT}, we define a torsion theory for any ring and some
related notions. We also introduce some examples of torsion
theories (Lambek, Goldie, classical, to name a few).

In Section \ref{TTforFVNAs}, we  study the torsion theories for a
finite von Neumann algebra $\A$. We introduce the torsion theory
$(\T,\bigP)$ for the algebra $\A$ (studied also in \cite{Lu1},
\cite{Lu2}, \cite{Re} for finitely generated modules) and show
that it is equal to both Lambek and Goldie torsion theories for
$\A$ (Proposition \ref{T=L=G}). This result will be needed in
Section \ref{EandK0}. The second torsion theory of interest is
$(\smallt, \p),$ where $\smallt$ is the class of modules which
vanishes when tensored with the algebra of affiliated operators
$\U$ of $\A.$ This torsion theory coincides with the classical
torsion theory for $\A$ (as an Ore ring). The torsion-free class
$\p$ consists of flat modules. The class $\smallt$ is the class of
cofinal-measurable modules studied in \cite{Lu4}, \cite{Lu5} and
\cite{Wegner}.

In Section \ref{EandK0}, we show that the module $\U\otimes_{\A}P$
is the injective envelope of a finitely generated projective
$\A$-module $P$ (Theorem \ref{envelope is tensor}). Using this, we
improve the known theorem on the isomorphism of $K_0(\A)$ and
$K_0(\U).$ Namely, it is known that the map
$[P]\mapsto[\U\otimes_{\A}P]$ for $P$ finitely generated
projective $\A$-module induces an isomorphism $K_0(\A)\rightarrow
K_0(\U).$ We show that the inverse of this isomorphism is the map
$[Q]\mapsto[Q\cap\A]$ (Theorem \ref{moja K 0 theorem}).

In Section \ref{capacity}, we study the class of finitely
generated $\A$-modules and capacity. Namely, for every
$\A$-module, there is a filtration
\[\underbrace{0\subseteq\smallt}_{\smallt
M}\underbrace{M\subseteq\T }_{\T\p M} \underbrace{M\subseteq
M.}_{\bigP M}\] We shall describe the three constitutive parts for
$M$ finitely generated in Proposition \ref{t-Tp-P}. The projective
quotient $\bigP M$ and the cofinal-measurable submodule $\smallt
M$ are faithfully measured by dimension and capacity,
respectively. We give a formula for calculating the capacity of a
module (Proposition \ref{my capacity}) that generalizes the one
from \cite{Lu4}.

In Section \ref{Induction}, we look at a subalgebra $\B$ of a
finite von Neumann algebra $\A.$ If $M$ is a $\B$-module, we
define the induced $\A$-module $i_*(M) = \A\otimes_{\B}M.$ In
\cite{Lu2}, it is shown that the dimension is invariant under the
induction of a module in the case of group von Neumann algebras.
We show that this holds for finite von Neumann algebras also. In
\cite{Wegner} it is shown that $c(M)\leq c(i_*(M))$ in the case of
a module over a group von Neumann algebra. We show that the
improved formula $c(M)=c(i_*(M))$ holds for any {\em finite} von
Neumann algebra (Theorem \ref{improved formula for c}).

\section{FINITE VON NEUMANN ALGEBRAS}
\label{FiniteVNA}

Let $H$ be a Hilbert space and $\B(H)$ be the algebra of bounded
operators on $H$. The space $\B(H)$ is equipped with five
different topologies: norm, strong, ultrastrong, weak and
ultraweak . The statements that a $*$-closed unital subalgebra
$\A$ of ${\B}(H)$ is closed in weak, strong, ultraweak and
ultrastrong topologies are equivalent (for details see \cite{D}).

\begin{defn} A {\em von Neumann algebra} $\A$ is a $*$-closed unital
subalgebra of ${\B}(H)$ which is closed with respect to weak
(equivalently strong, ultraweak, ultrastrong) operator topology.
\end{defn}

A $*$-closed unital subalgebra $\A$ of ${\B}(H)$ is a von Neumann
algebra if and only if $\A= \A''$ where $\A'$ is the commutant of
$\A.$ The proof can be found in \cite{D}.

\begin{defn} A von Neumann algebra $\A$ is {\em finite} if
there is a $\Cset$-linear function $\tr_{\A}: \A\rightarrow\Cset$
such that
\begin{enumerate}
\item $\tr_{\A}(ab) =\tr_{\A}(ba).$
\item $\tr_{\A}(a^*a)\geq 0.$ If a $\Cset$-linear function on $\A$
satisfies 1. and 2., we call it a {\em trace}.
\item $\tr_{\A}$ is {\em normal}: it is continuous with respect to
ultraweak topology.
\item $\tr_{\A}$ is {\em faithful}: $\tr_{\A}(a)=0$ for some $a\geq 0$
(i.e. $a = bb^*$ for some $b\in\A$) implies $a=0$.
\end{enumerate}
 \label{finite vna}
\end{defn}

Trace is not unique. A trace function extends to matrices over
$\A$ in a natural way: the trace of a matrix is the sum of the
traces of the elements on the main diagonal.

\begin{exmp}
Let $G$ be a (discrete) group. The group ring $\Cset G$ is a
pre-Hilbert space with an inner product: $ \langle \;\sum_{g\in G}
a_g g, \sum_{h\in G} b_h h\;\rangle = \sum_{g\in
G}a_g\overline{b_g}.$

Let $\lge$ be the Hilbert space completion of $\Cset G.$ Then
$\lge$ is the set of square summable complex valued functions over
the group $G$.

{\em The group von Neumann algebra} $\vng$ is the space of
$G$-equivariant bounded operators from $\lge$ to itself: \[\vng =
\{\; f\in {\B}(\lge)\;|\;f(gx)=gf(x)\mbox{ for all }g\in G\mbox{
and }x\in\lge\;\}.\]

$\Cset G$ embeds into ${\B}(\lge)$ by right regular
representations. $\vng$ is a von Neumann algebra for $H = \lge$
since it is the weak closure of $\Cset G$ in ${\B}(\lge)$ so it is
a $*$-closed subalgebra of ${\B}(\lge)$ which is weakly closed
(see Example 9.7 in \cite{Lu5} for details). $\vng$ is finite as a
von Neumann algebra since it has a normal, faithful trace
$\tr_{\A}(f) = \langle f(1), 1 \rangle.$
\end{exmp}

One of the reasons a finite von Neumann algebra is attractive is
that the trace provides us with a way of defining a convenient
notion of the dimension of any module.

\begin{defn}
If $P$ is a finitely generated projective $\A$-module, there exist
$n$ and $f:\A ^n\rightarrow\A^n$ such that $f=f^2=f^*$ and the
image of $f$ is $P.$ Then, the {\em dimension} of $P$ is \[
\dim_{\A}(P)=\tr_{\A}(f).\]

If $M$ is any $\A$-module, the {\em dimension} $\dim_{\A}(M)$ is
defined as \[\dim_{\A}(M)=\sup \{ \dim_{\A}(P)\; |\; P \mbox{
finitely generated projective submodule of }M\}.\]
\label{dimension of any vng module}
\end{defn}

In the first part of the definition, the map $f^*$ is defined by
transposing and applying $*$ to each entry of the matrix
corresponding to $f.$

The dimension of a finitely generated projective $\A$-module is a
nonnegative real number, while the dimension of any $\A$-module is
in $[0, \infty].$ The dimension has the following properties.
\begin{prop}
\begin{enumerate}
\item If $\;0\rightarrow M_0\rightarrow M_1\rightarrow M_2\rightarrow
0$ is a short exact sequence of $\A$-modules, then $
\dim_{\A}(M_1)= \dim_{\A}(M_0)+\dim_{\A}(M_2).$
\item If $M = \bigoplus_{i\in I} M_i,$ then $\dim_{\A}(M) =
\sum_{i \in I}\dim_{\A}(M_i).$
\item If $M = \bigcup_{i \in I}M_i$ is a directed union of
submodules, then $\dim_{\A}(M) = \sup\{\;\dim_{\A}(M_i)\; |\; i\in
I\;\}.$
\item If $M = \dirlim_{i\in I} M_i$ is a direct limit of a directed
system, then $\dim_{\A}(M) \leq \lim\inf\{\dim_{\A}(M_i)\;|\;i\in
I\}.$
\item If $M$ is a finitely generated projective module, then
$\;\dim_{\A}(M)=0$ iff $M=0.$
\end{enumerate}
\label{propofdim}
\end{prop}

The proof of this proposition can be found in \cite{Lu2}.

Besides using the above approach to derive the notion of the
dimension of an $\A$-module, we can use a more operator theory
oriented approach.

A finite von Neumann algebra is a pre-Hilbert space with inner
product $\langle a, b \rangle = \tr_{\A}(ab^*).$ Let $\la$ denote
the Hilbert space completion of $\A.$

Clearly, this is the analogue of $\lge$ for $\A$ in case when $\A$
is the group von Neumann algebra $\vng$ of some group $G$. This is
because $\lge,$ as defined before, is isomorphic to $l^2(\vng),$
as defined above. These two spaces are isomorphic since they are
both Hilbert space completions of $\vng$ (see section 9.1.4 in
\cite{Lu5} for details).

A finite von Neumann algebra $\A$ can be identified with the set
of $\A$-equivariant bounded operators on $\la,$ ${\B}(\la)^{\A},$
using the right regular representations. This justifies our
definition of $\vng$ as $G$-equivariant operators in ${\B}(\lge)$
since $ {\B}(l^2(\vng))^{\vng}= {\B}(\lge)^{\vng} = {\B}(\lge)^{G}
= \vng.$

In \cite{Lu1} the following theorem is proved.

\begin{thm} There is an equivalence of categories
\[\nu :\{ \mbox{fin. gen. proj. }\A\mbox{-mod.} \}\rightarrow\{
\mbox{fin. gen. Hilbert }\A\mbox{-mod.}\}\] \label{equivalence}
\end{thm}

Here, a {\em finitely generated Hilbert $\A$-module} is a Hilbert
space $V$ with a left representation $\A\rightarrow \B(V)$ and
such that there is a nonnegative integer $n$ and a projection
$p:(\la)^n$ $\rightarrow$ $(\la)^n$ whose image is isometrically
$\A$-isomorphic to $V.$ Such projection $p$ can be viewed as an
$n\times n $ $\A$-matrix (when identifying ${\B}(\la)^{\A}$ and
$\A$). The dimension of such $V$ is defined as $\tr_{\A}(p).$

Using the above equivalence of the categories, we can define the
dimension of a finitely generated projective $\A$-module $P$ via
the dimension of a finitely generated Hilbert $\A$-module $\nu(P)$
as \[\dim_{\A}(P) = \dim_{\A}(\nu(P)) \] This definition agrees
with the first part of Definition \ref{dimension of any vng
module}. The dimension defined in this way for finitely generated
projective modules extends to all modules in the same way as in
Definition \ref{dimension of any vng module}.

Theorem \ref{equivalence} allows us to choose between a more
algebraic and a more operator theory oriented approach. This is
just one example of the accord between algebra and operator theory
related to finite von Neumann algebras. There will be other
examples of this phenomenon later on.

Let us turn to some ring-theoretic properties of a finite von
Neumann algebra $\A$. As a ring, $\A$ is {\em semihereditary}
(i.e., every finitely generated submodule of a projective module
is projective or, equivalently, every finitely generated ideal is
projective). This follows from two facts. First, every von Neumann
algebra is a Baer $*$-ring and, hence, a Rickart $C^*$-algebra
(see Chapter 1.4 in \cite{Be2}). Second, a $C^*$-algebra is
semihereditary as a ring if and only if it is Rickart (see
Corollary 3.7 in \cite{AG}).

Alternative proof of the fact that $\A$ is semihereditary uses the
equivalence from Theorem \ref{equivalence}. It can be found in
\cite{Lu1}.

$\A$ is also a (left and right) {\em nonsingular} ring. Recall
that a ring $R$ is left nonsingular if, for every $x\in R,$
ann$_l(x) = \{r\in R\;|\; rx = 0\}$ intersects every left ideal
nontrivially if and only if $x = 0.$ The right nonsingular ring is
defined analogously. $\A$ is nonsingular as a Rickart ring (see
7.6 (8) and 7.48 in \cite{Lam}). Alternatively, $\A$ is
nonsingular since it is a $*$-ring with involution such that
$x^*x=0$ implies $x=0$ (see 7.9 in \cite{Lam}).

\subsection{The Algebra of Affiliated Operators}

\begin{defn} Let $a$ be a linear map
$a: \mathrm{dom}\; a$ $\rightarrow$ $\la$, $\mathrm{dom}\; a
\subseteq\la.$ We say that $a$ is {\em affiliated to $\A$} if
\begin{itemize}
\item[i)] $a$ is densely defined (the domain $\mathrm{dom}\;
 a$ is a dense subset of $\la);$
\item[ii)] $a$ is closed (the graph of $a$ is closed in
$\la \oplus \la);$
\item[iii)]  $ba = ab$ for every $b$ in the commutant of $\A.$
\end{itemize}
Let $\U=\U(\A)$ denote the {\em algebra of operators affiliated
to} $\A$.
\end{defn}

$\U(\A)$ is an $*$-algebra with $\A$ as a $*$-subalgebra.

\begin{prop}
Let $\A$ be a finite von Neumann algebra and $\U = \U(\A)$ its
algebra of affiliated operators.
\begin{enumerate}
\item $\A$ is an Ore ring.
\item $\U$ is equal to the classical ring of quotients
$Q_{\mathrm{cl}}(\A)$ of $\A.$
\item $\U$ is a von Neumann regular (fin. gen. submodule of
fin. gen. projective module is a direct summand), left and right
self-injective ring.
\end{enumerate}
\label{ug is clasical and maximal}
\end{prop}
The proof of 1. and 2. can be found in \cite{Re}. The proof of 3.
can be found in \cite{Be1}.

From this proposition and the fact that $\A$ is a semihereditary
ring, it follows that the algebra of affiliated operators $\U$ is
both the maximal $Q_{\mathrm{max}}(\A)$ and the classical
$Q_{\mathrm{cl}}(\A)$ ring of quotients of $\A$ as well as the
injective envelope $E(\A)$ of $\A$ (minimal injective module
containing $\A$). Thus, the algebra $\U$ can be defined both by
using purely algebraic terms (ring of quotient, injective
envelope) and by using just the operator theory terms (affiliated
operators).

The ring $\U$ has many nice properties that $\A$ is missing: it is
von Neumann regular and self-injective; and it keeps all the
properties that $\A$ has: it is semihereditary and nonsingular.

\section{TORSION THEORIES}
\label{TT}

The ring $\A$ is very handy to work with because it has many
PID-like features. Every finitely generated module over a
principal ideal domain (PID) is the direct sum of its torsion and
torsion-free part. Our ring $\A$ has the similar property.
However, for the ring $\A$ it turns out that there exists more
than just one natural definition of a torsion element. To study
the different ways to define the torsion and torsion-free part of
an $\A$-module, we first introduce the general framework in which
we shall be working
--- the torsion theory.

\begin{defn} Let $R$ be any ring. A {\em torsion theory}
for $R$ is a pair $\tau = (\te, \ef)$ of classes of $R$-modules
such that
\begin{itemize}
\item[i)] $\homo_R(T,F)=0,$ for all $T \in \te$ and $F \in \ef.$
\item[ii)] $\te$ and $\ef$ are maximal classes having the property
$i).$
\end{itemize}
\end{defn}
The modules in $\te$ are called {\em $\tau$-torsion modules} (or
torsion modules for $\tau$) and the modules in $\ef$ are called
{\em $\tau$-torsion-free modules} (or torsion-free modules for
$\tau$).

If $\tau_1 = (\te_1, \ef_1)$ and $\tau_2 = (\te_2, \ef_2)$ are two
torsion theories, we say that $\tau_1$ is {\em smaller} than
$\tau_2$ \[\tau_1\leq\tau_2\mbox{ iff }\te_1\subseteq\te_2\mbox{
iff }\ef_1\supseteq\ef_2.\]

If $\ce$ is a class of $R$-modules, then torsion theory {\em
generated} by $\ce$ is the smallest torsion theory $(\te, \ef)$
such that $\ce\subseteq\te.$

The torsion theory {\em cogenerated} by $\ce$ is the largest
torsion theory $(\te, \ef)$ such that $\ce\subseteq\ef.$

\begin{prop}
\begin{enumerate}
\item If $(\te, \ef)$ is a torsion theory,
\begin{itemize}
\item[i)] the class $\te$ is closed under quotients, direct sums and
extensions;
\item[ii)]  the class $\ef$ is
closed under taking submodules, direct products and extensions.
\end{itemize}
\item  If $\ce$ is a class of $R$-modules closed under quotients,
direct sums and extensions, then it is a torsion class for a
torsion theory $(\ce, \ef)$ where $\ef =
\{\;F\;|\;\homo_R(C,F)=0,\mbox{ for all }C\in\ce\;\}.$

Dually, if $\ce$ is a class of $R$-modules closed under
submodules, direct products and extensions, then it is a
torsion-free class for a torsion theory $(\te, \ce)$ where $\te =
\{\;T\;|\;\homo_R(T,C)=0,\mbox{ for all }C\in\ce\;\}.$
\item Two classes of $R$-modules $\te$ and $\ef$ constitute a
torsion theory iff
\begin{itemize}
\item[i)] $\te\cap\ef = \{0\},$
\item[ii)] $\te$ is closed under quotients,
\item[iii)] $\ef$ is closed under submodules and
\item[iv)] For every module $M$ there exists submodule $N$ such that
$N\in\te$ and $M/N\in\ef.$
\end{itemize}
\end{enumerate}
\end{prop}
The proof of this proposition is straightforward by the definition
of a torsion theory. The details can be found in \cite{Bland}. In
iv) of part (3) take $N$ to be the submodule of $M$ generated by
the union of all torsion submodules of $M$.

From this proposition it follows that every module $M$ has the
largest submodule which belongs to $\te$. We call it the {\em
torsion submodule} of $M$ and denote it with $\te M$. The quotient
$M/\te M$ is called the {\em torsion-free quotient} and we denote
it $\ef M.$

We say that a torsion theory $\tau = (\te, \ef)$ is {\em
hereditary} if the class $\te$ is closed under taking submodules.
A torsion theory is hereditary if and only if the torsion-free
class is closed under formation of injective envelopes. Also, a
torsion theory cogenerated by a class of injective modules is
hereditary and, conversely, every hereditary torsion theory is
cogenerated by a class of injective modules. The proof of these
facts is straightforward. The details can be found in
\cite{Golan}.

Some authors (e.g. \cite{Golan}) consider just hereditary torsion
theories and call a torsion theory what we here call a hereditary
torsion theory.

The notion of the closure of a submodule in a module is another
natural notion that can be related to a torsion theory.
\begin{defn}
If $M$ is an $R$-module and $K$ a submodule of $M,$ let us define
the {\em closure} $\cl_{\te}^M(K)$ of $K$ in $M$ with respect to
the torsion theory $(\te, \ef)$ by \[\cl_{\te}^M(K) =
\pi^{-1}(\te(M/K))\mbox{ where } \pi\mbox{ is the natural
projection }M\twoheadrightarrow M/K.\]
\end{defn}

If it is clear in which module we are closing the submodule $K,$
we suppress the superscript $M$ from $\cl_{\te}^M(K)$ and write
just $\cl_{\te}(K)$. If $K$ is equal to its closure in $M,$ we say
that $K$ is {\em closed} submodule of $M$.

The closure has the following properties.
\begin{prop}
Let $(\te, \ef)$ be a torsion theory on $R$, let $M$ and $N$ be
$R$-modules and $K$ and $L$ submodules of $M$. Then
\begin{enumerate}
\item $\te M = \cl_{\te}(0).$
\item $\te(M/K) = \cl_{\te}(K)/K$ and
$\ef(M/K) \cong M/\cl_{\te}(K).$
\item If $K \subset L,$ then $\cl_{\te}(K)\subseteq\cl_{\te}(L).$
\item $K\subset\cl_{\te}(K)$ and
$\cl_{\te}(\cl_{\te}(K)) = \cl_{\te}(K).$
\item $\cl_{\te}(K)$ is the smallest closed submodule of
$M$ containing $K.$
\item If $(\te, \ef)$ is hereditary, then
$\cl_{\te}^K(K\cap L) = K \cap\cl^M_{\te}(L).$ If $(\te, \ef)$ is
not hereditary, just $\subseteq$ holds in general.
\item If $(\te_1, \ef_1)$ and $(\te_2, \ef_2)$ are two torsion
theories, then \[(\te_1, \ef_1)\leq(\te_2, \ef_2)\;\;\; \mbox{ iff
}\;\;\; \cl_{\te_1}(K)\subseteq\cl_{\te_2}(K)\;\mbox{ for all
}K.\]
\end{enumerate}
\label{properties of closure}
\end{prop}

The proof follows directly from the definition of the closure.

\subsection{Examples}
\label{Examples}

\begin{enumerate}

\item {\em The trivial torsion theory} for $R$ is the torsion
theory $(0, \mathrm{Mod}_R),$ where Mod$_R$ is the class of all
$R$-modules.

\item {\em The improper torsion theory} for $R$ is the torsion
theory $(\mathrm{Mod}_R, 0).$

\item The torsion theory cogenerated by the injective envelope
$E(R)$ of $R$ is called the {\em Lambek torsion theory}. We denote
it $\tau_L.$ Since it is cogenerated by an injective module, it is
hereditary.

If the ring $R$ is torsion-free in a torsion theory $\tau$, we say
that $\tau$ is {\em faithful}. It is easy to see that the Lambek
torsion theory is faithful. Moreover, it is the largest hereditary
faithful torsion theory. For more details, see \cite{Tele_Teza}.

\item The class of nonsingular modules over a ring $R$ is closed
under submodules, extensions, products and injective envelopes.
Thus, the class of all nonsingular modules is a torsion-free class
of a hereditary torsion theory. This theory is called the {\em
Goldie torsion theory} $\tau_G.$

The Lambek theory is smaller than the Goldie theory. This is the
case since the Goldie theory is larger than any hereditary torsion
theory (see \cite{Bland}). Moreover, the Lambek and Goldie
theories coincide if and only if $R$ is a nonsingular ring (i.e.
$\tau_G$ is faithful). For more details, see \cite{Tele_Teza}.

A finite von Neumann algebra is a nonsingular ring so its Lambek
and Goldie torsion theories coincide.

\item
If $R$ is an Ore ring with the set of regular elements $T$ (i.e.,
$Tr \cap Rt \neq 0,$ for every $t \in T$ and $r\in R$), we can
define a hereditary torsion theory by the condition that an
$R$-module $M$ is a torsion module iff for every $m\in M$, there
is a nonzero $t\in T$ such that $tm =0.$ This torsion theory is
called the {\em classical torsion theory of an Ore ring}. This
theory is faithful and so it is contained in the Lambek torsion
theory.

Note the following lemma that we will use in the sequel.
\begin{lem} If $t$ is a regular element of an Ore ring $R$ and $r\in
R,$ then $t r t^{-1}$ is in $R.$ \label{regular_elements_lemma}
\end{lem}
Note that $r t^{-1}$ is defined as an element of the classical
ring of quotients $Q_{\mathrm{cl}}(R).$
\begin{proof} Since $t$ is regular, $tR=R=Rt.$ Thus, $tr=st$ for
some $s\in R.$ Hence, $t r t^{-1}= st t^{-1} = s\in R.$
\end{proof}

\item Let $R$ be a subring of a ring $S$. Let us look at a
collection $\te$ of $R$-modules $M$ such that $S\otimes_R M = 0.$
This collection is closed under quotients, extensions and direct
sums. Moreover, if $S$ is flat as an $R$-module, then $\te$ is
closed under submodules and, hence, defines a hereditary torsion
theory. In this case, denote this torsion theory by $\tau_S.$

From the definition of $\tau_S$ it follows that
\begin{itemize}
\item[1.] The torsion submodule of $M$ in $\tau_S$ is the kernel
of the natural map $M\rightarrow S \otimes_R M,$ i.e.
$\tor^R_1(S/R, M).$
\item[2.] All flat modules are $\tau_S$-torsion-free.
\end{itemize}
By 2., $\tau_S$ is faithful. Thus, $\tau_S$ is contained in the
Lambek torsion theory.

If a ring $R$ is Ore, then the classical ring of quotients
$Q_{\mathrm{cl}}^l(R)$ is a flat $R$-module and the set $\{\;m\in
M\;|\; rm = 0, $ for some nonzero-divisor $r\in R\;\}$ is equal to
the kernel $\ker (M\rightarrow Q_{\mathrm{cl}}^l(R) \otimes_R M).$
Hence the torsion theory $\tau_{Q_{\mathrm{cl}}^l(R)}$ coincides
with the classical torsion theory of $R$ in this case.

Since $\A$ is Ore and $\U=Q_{\mathrm{cl}}(\A),$ $\U$ is a flat
$\A$-module and $\tau_{\U}$ is the classical torsion theory of
$\A.$

\item All the torsion theories we introduced so far are hereditary.
Let us introduce a torsion theory that is not necessarily
hereditary. Let $(\bnd, \unb)$ be the torsion theory cogenerated
by the ring $R$. This is the largest torsion theory in which $R$
is torsion-free. We call a module in $\bnd$ a {\em bounded module}
and a module in $\unb$ an {\em unbounded module}.

The Lambek torsion theory $\tau_L$ is contained in the torsion
theory $(\bnd, \unb)$ because $R$ is $\tau_L$-torsion-free. There
is another interesting relation between the Lambek and $(\bnd,
\unb)$ torsion theory. Namely,
\begin{center}
$M$ is a Lambek torsion module if and only if every submodule of
$M$ is bounded.
\end{center}
This is a direct corollary of the fact that $\homo_R(M, E(R)) = 0$
if and only if $\homo_R(N, R)=0,$ for all submodules $N$ of $M,$
which is an exercise in \cite{Cohn1}. Also, it is easy to show
that $(\bnd, \unb)$ is equal to the Lambek torsion theory if and
only if $(\bnd, \unb)$ is hereditary.
\end{enumerate}

To summarize, for any ring $R$ we have the following relationship
for the torsion theories:
\begin{center} Trivial $\leq$ Lambek $\leq$ Goldie
$\leq$ $(\bnd, \unb)$ $\leq$ Improper.
\end{center}
If $R$ is an Ore nonsingular ring, then
\begin{center}
Trivial $\leq$ Classical = $\tau_{Q_{\mathrm{cl}}(R)}$ $\leq$
Lambek $=$ Goldie $\leq$ $(\bnd, \unb)$ $\leq$ Improper.
\end{center}
The last is the situation for our finite von Neumann algebra $\A$.
We shall examine the situation further in the next section.

\section{TORSION THEORIES FOR FINITE VON NEUMANN ALGEBRAS}
\label{TTforFVNAs}

In this section, we shall introduce some theories for group von
Neumann algebras and identify some of them with the torsion
theories from previous section.

\begin{enumerate}
\item The dimension of an $\A$-module
enables us to define a torsion theory. For an $\A$-module $M$
define $\T M$ as the submodule generated by all submodules of $M$
of the dimension equal to zero. It is zero-dimensional by property
3. (Proposition \ref{propofdim}). So, $\T M$ is the largest
submodule of $M$ of zero dimension.

Let us denote the quotient $M/ \T M$ by $\bigP M$.

Proposition \ref{propofdim} gives us that the class $ \T = \{M \in
\mathrm{Mod}_{\A} | M = \T M\}$ is closed under submodules,
quotients, extensions and direct sums. Thus, $\T$ defines a
hereditary torsion theory with torsion-free class equal to $\bigP
= \{M \in \mathrm{Mod}_{\A} | M = \bigP M\}.$

From the definition of this torsion theory it follows that
$\cl_{\T}(K)$ is the largest submodule of $M$ with the same
dimension as $K$ for every submodule $K$ of an module $M$. Also,
since $\A$ is semihereditary and a nontrivial finitely generated
projective module has nontrivial dimension, $\A$ is in $\bigP$ and
so the torsion theory $(\T, \bigP)$ is faithful.

\item The second torsion theory of interest is $(\bnd, \unb),$ the
largest torsion theory in which the ring is torsion-free. Since
$\A$ is torsion-free in $(\T, \bigP),$ we have that $(\T, \bigP)
\leq (\bnd, \unb).$

\item Let $(\smallt, \p)$ denote the classical torsion theory of
$\A.$ Since $\U=Q_{\mathrm{cl}}(\A),$ \[\smallt M = \ker (M
\rightarrow \U \otimes_{\A}M)= \tor_1^{\A}(\U/\A, M)\] for any
$\A$-module $M$ (see Examples (5) and (6) in Subsection
\ref{Examples}). We denote the torsion-free quotient $M/ \smallt
M$ by $\p M.$ From Example (6), it follows that all flat modules
are torsion-free. In \cite{Turnidge}, the torsion theory from
example (6) is studied. Turnidge showed in \cite{Turnidge} that
all torsion-free modules are flat if the following conditions
hold:
\begin{itemize}
\item[-] The ring $R$ is semihereditary;
\item[-] The ring $Q$ is von Neumann regular;
\item[-] $Q$ is flat as an $R$-module.
\end{itemize}

The finite von Neumann algebra $\A$ is semihereditary, $\U$ is von
Neumann regular, and $\A$-flat. Thus for an $\A$-module $M$ the
following is true
\begin{center}
$M$ is flat if and only if $M$ is in $\p.$
\end{center}

The class of flat modules of a semihereditary ring is closed under
submodules, extensions and direct product and, hence, is a
torsion-free class of a torsion theory. Turnidge's theorem states
that this torsion theory is exactly the classical torsion theory
$(\smallt, \p).$

It turns out that the torsion class $\smallt$ also demonstrates
the accord between the algebra and operator theory ingrained in
$\A.$ Namely, the class $\smallt$ (defined as above using purely
algebraic notions) coincide with the class of cofinal-measurable
modules defined using the dimension function and hence operator
theory. We say that an $\A$-module $M$ is {\em measurable} if it
is a quotient of a finitely presented module of dimension zero.
$M$ is {\em cofinal-measurable} if each finitely generated
submodule is measurable. The class $\smallt$ is the class of
cofinal-measurable modules. For proof of this fact, see \cite{Lu5}
(proof is given for a group von Neumann algebra but it holds for
any finite von Neumann algebra).
\end{enumerate}

Let us now compare the defined torsion theories.

Example 8.35 in \cite{Lu5} shows that $\T$ is different than
$\bnd$ in general. Still, the torsion theories $(\T, \bigP)$ and
$(\bnd, \unb)$ coincide on finitely generated modules as the
following proposition shows.

\begin{prop}
Let $M$ be a finitely generated $\A$-module and $K$ a submodule of
$M$. Then
\begin{itemize}
\item[i)] $\dim_{\A}(K) = \dim_{\A}(\cl_{\bnd}(K)).$
\item[ii)] $\cl_{\bnd}(K)$ is a direct summand in $M$ and
$M/\cl_{\bnd}(K)$ is finitely generated projective module.
\item[iii)] $\cl_{\T}(K)=\cl_{\bnd}(K).$ In particular,
$\T M = \bnd M.$
\item[iv)] $M = \T M\oplus\bigP M = \bnd M \oplus \unb M,$
\end{itemize}
\label{rezultat iz Lucka dim}
\end{prop}

The proof of i) and ii) can be found in \cite{Lu2}. The idea of
the proof is to first show i) and ii) for a special case when $M$
is projective. In this case, the proposition is proven using the
equivalence of the category of finitely generated projective
$\A$-modules and the finitely generated Hilbert $\A$-modules
(Theorem \ref{equivalence}). Then the general case is proven.

To prove part iii), note that from part i) it follows that
$\cl_{\bnd}(K)\subseteq\cl_{\T}(K)$ because $\cl_{\T}(K)$ is the
largest submodule of $M$ containing $K$ with the same dimension as
$K$. But since $\T\subseteq\bnd,$ the converse holds as well.
Thus, $\cl_{\T}(K)=\cl_{\bnd}(K).$ Taking $K=0$ gives us $\T M
=\bnd M.$

Part iv) follows from ii) and iii).

This proposition gives us that every finitely generated module in
$\bigP$ is projective. This gives us a nice characterization of
any module in $\bigP.$ Namely,
\begin{center} an $\A$-module $M$ is a $\bigP$-module iff
every finitely generated submodule of $M$ is projective.
\end{center}
Thus, a $\bigP$ module is a directed union of finitely generated
projective modules.

\begin{prop}
For the ring $\A,$
\begin{center} $(\T, \bigP) =$ Lambek torsion theory $=$ Goldie
torsion theory.
\end{center}
\label{T=L=G}
\end{prop}
\begin{proof}
The Lambek torsion theory $\tau_L$ is the same as the Goldie
torsion theory $\tau_G$ because $\A$ is a nonsingular ring. Since
$\tau_L$ is the largest hereditary torsion theory in which the
ring is torsion-free and $\A$ is torsion-free in $(\T, \bigP),$ we
have that $(\T, \bigP)\leq$ $\tau_L$ = $\tau_G.$

To prove the first equality, we shall prove that every Lambek
torsion module $M$ has dimension zero. Recall that $M$ is Lambek
torsion module iff all submodules of $M$ are bounded. This means
that all finitely generated submodules of $M$ are in $\T$ (a
finitely generated module is in $\bnd$ iff it is in $\T$ by
Proposition \ref{rezultat iz Lucka dim}). The dimension of $M$ is
equal to the supremum of the dimensions of finitely generated
submodules of $M$ by Proposition \ref{propofdim}. But that
supremum is 0, so $M$ is in $\T.$
\end{proof}

This proposition is another example of the harmony between algebra
and the operator theory in a finite von Neumann algebra $\A$. The
proposition asserts that the theory $(\T, \bigP)$ (defined using
the dimension i.e. the operator theory) is the same theory as the
Goldie or Lambek theories, the theories defined via purely
algebraic notions.

It is also interesting that this proposition shows that the
torsion theory $(\T, \bigP)$, defined via a normal and faithful
trace $\tr_{\A},$ is not dependent on the choice of such trace
since $(\T, \bigP)$ coincides with Lambek and Goldie theories.

Let us compare the theory $(\smallt,\p)$ with the other torsion
theories of $\A.$ Since $\A$ is flat as $\A$-module, the ring $\A$
is torsion-free in $(\smallt, \p)$. Hence, this torsion theory is
contained in $\tau_L$ (recall that the Lambek torsion theory is
the largest hereditary theory in which the ring is torsion-free).
But $\tau_L$ is the same as $(\T, \bigP),$ and so we have
$(\smallt, \p)\leq (\T, \bigP).$ The examples that $\smallt M
\subsetneq \T M$ can be found even for $M$ finitely generated
(Example 8.34 in \cite{Lu5}). However, the classes $\T$ and
$\smallt$ coincide when restricted on the class of finitely
presented $\A$-modules (see Lemma 8.33 in \cite{Lu5}).

The theory $(\smallt,\p)$ can be nontrivial by Example 2.9 in
\cite{Lu5}.

For any nontrivial finite von Neumann algebra $\A$, the theory
$(\bnd, \unb)$ is not improper since $\A$ is a module in $\unb.$

To summarize, various torsion theories for $\A$ are ordered as
follows:
\begin{center}
Trivial $\leq$ Classical = $\tau_{\U}$ = $(\smallt,\p)$ $\leq$
$\tau_{L}$ $=$ $\tau_{G}$ = $(\T, \bigP)$ $\leq$ $(\bnd, \unb)$
$\leq$ Improper
\end{center}
where all of the above inequalities can be strict.

The following proposition further explores the relations between
the torsion theories for $\A.$

\begin{prop}
\begin{enumerate}
\item $\T\smallt=\smallt\T=\smallt,$ $\smallt\bigP=\bigP\smallt=0,$
and $\p\bigP=\bigP;$
\item $\cl_{\T}(\smallt M) = \T M$ for every $\A$-module $M$;
\item $\bigP\p \cong \bigP;$
\item $\T\p \cong \p\T.$
\end{enumerate}
\label{odnosi Tova i Pova}
\end{prop}
\begin{proof}
The equations in (1) are direct consequences of the fact that
$\smallt\subseteq\T$ and that the torsion and torsion-free classes
intersect trivially.

(2) Since $\smallt M$ has dimension zero, $\cl_{\T}(\smallt M)$
has dimension zero as well. So, $\cl_{\T}(\smallt M)\subseteq \T
M.$ The other inclusion follows since $\T M = \cl_{\T}(0)\subseteq
\cl_{\T}(\smallt M).$

(3) $\bigP \p M  =  \bigP(M/\smallt M) \cong M / \cl_{\T}(\smallt
M)$ by Proposition \ref{properties of closure}.
$M/\cl_{\T}(\smallt M)= M/\T M = \bigP M$ by (2) above.

(4) We shall show that both $\T\p M$ and $\p\T M$ are isomorphic
to the quotient $\T M/\smallt M$. First, $\T\p M = \T(M/\smallt M)
= \cl_{\T}(\smallt M)/\smallt M = \T M/\smallt M.$ We obtain the
middle equality by Proposition \ref{properties of closure} and the
last one by (2) above.

$\p\T M$ is isomorphic to $\T M / \cl_{\smallt}^{\T M}(\smallt M)$
by Proposition \ref{properties of closure}. But the closure of
$\smallt M$ with respect to $(\smallt, \p)$ is the same both in
$M$ and in $\T M$ since $\cl_{\smallt}^{\T M}(\smallt M) = \T M
\cap \smallt M = \smallt M = \cl_{\smallt}^M(\smallt M).$ Thus,
$\p\T M \cong \T M / \cl_{\smallt}^{\T M}(\smallt M) = \T
M/\smallt M.$

\end{proof}

This proposition gives us that for every module $M,$ we have a
filtration: \[\underbrace{0\subseteq\smallt}_{\smallt
M}\underbrace{M\subseteq\T }_{\T\p M} \underbrace{M\subseteq
M.}_{\bigP M}\]

Thus, every module is built up of three building blocks:
\begin{itemize}
\item[1.] a cofinal-measurable part $\smallt M$,
\item[2.] a flat, zero-dimensional part $\p\T M = \T\p M,$
\item[3.] a $\bigP$-part $\bigP M,$
(directed union of finitely generated projective modules;
projective if finitely generated).
\end{itemize}

If $M$ is finitely presented, $\p\T M=0,$ (since $\smallt M = \T
M$) so there are just two parts: $\T M = \smallt M$ and $\bigP M,$
and they are direct summands of $M$.

For $M$ finitely generated, $\p\T M$ does not have to vanish
(Example 8.34 in \cite{Lu5}) but the finitely generated quotient
$\p M$ splits as the direct sum of $\T\p M$ and $\bigP\p M = \bigP
M$ and thus we have a short exact sequence $ 0\rightarrow\smallt
M\rightarrow M\rightarrow \T\p M\oplus\bigP M \rightarrow 0.$

\section{INJECTIVE ENVELOPES AND $K_0$-THEOREM}
\label{EandK0}

In this section, we shall obtain some results on the injective
envelopes of $\A$-modules and show that the injective envelope of
a finitely generated projective module $P$ is $\U\otimes_{\A}P.$
Using that, we shall acquire some further results on the
isomorphism on $K_0$ of $\A$ and $\U.$ Namely, Handelman proved
(Lemma 3.1 in \cite{Handel}) that for every finite Rickart
$C^*$-algebra $\A$ such that every matrix algebra over $\A$ is
also Rickart, the inclusion of $\A$ into a certain regular ring
$R$ with the same lattice of projections as $\A$ induces an
isomorphism $\mu: K_0(\A)\rightarrow K_0(R).$

By Theorem 3.4 in \cite{AG}, a matrix algebra over a Rickart
$C^*$-algebra is a Rickart $C^*$-algebra. Thus, $K_0(\A)$ is
isomorphic to $K_0(R)$ for every finite Rickart $C^*$-algebra. If
$\A$ is a finite von Neumann algebra, the ring $R$ can be
identified with the maximal ring of quotients
$Q_{\mathrm{max}}(\A)$ (e.g. \cite{Be1} and \cite{Be2}). This
gives us that the inclusion of a finite von Neumann algebra $\A$
in its algebra of affiliated operators $\U$ induces the
isomorphism \[\mu: K_0(\A)\rightarrow K_0(\U).\] Here, we shall
obtain the explicit description of the map Proj($\U$)
$\rightarrow$ Proj($\A$) that induces the inverse of the
isomorphism $\mu$ on $K_0$ of $\A$ and $\U.$

\subsection{Preliminaries}

Let $R$ be any ring. A submodule $K$ of an $R$-module $M$ is an
{\em essential submodule} of $M$ ($K\subseteq_e M$) if $K\cap
L\neq 0,$ for every nonzero submodule $L$ of $M.$ If $K\subseteq_e
M$, $M$ is an {\em essential extension} of $K$. $M$ is a maximal
essential extension of $K$ if no module strictly containing $M$ is
essential extension of $K$. Besides being defined as the minimal
injective module containing $M$, the injective envelope $E(M)$ can
be defined as a unique (up to isomorphism) maximal essential
extension of $M$. Hence, $M\subseteq_e E(M).$

A submodule $K$ of $M$, is a {\em complement} in $M$
($K\subseteq_c M$) if there exists a submodule $L$ of $M$ such $K$
is a maximal submodule of $M$ with the property that $K\cap L =
0.$ We shall use the following proposition from \cite{Lam}
(Proposition 6.32 in \cite{Lam}).

\begin{prop}
If $M$ is an $R$-module and $K$ a submodule of $M$, then the
following are equivalent:
\begin{itemize}
\item[a)] $K\subseteq_c M$
\item[b)] $K$ does not have any proper essential extensions in $M$.
\item[c)] $K$ is the intersection of M with a direct summand of $E(M).$
\end{itemize}
Moreover, if $L$ is a direct summand of $E(M)$ then $K=L\cap M$
satisfies a)-c). \label{Lamova prop o closed submodules}
\end{prop}

From the proof of this proposition it follows that if
$K\subseteq_c M,$ then the direct summand from part c) of the
above proposition is $E(K)$ and $K = E(K)\cap M.$

If $R$ is a nonsingular ring, we can describe the closure of a
submodule of nonsingular module with respect to the Goldie torsion
theory via the notion of an essential extension. Namely, the
following proposition holds.

\begin{prop}
Let $R$ be a nonsingular ring, $M$ an $R$-module and $K$ a
submodule of $M$. Then, the Goldie closure of $K$ in $M$ is
complemented in $M$. If $M$ is nonsingular, then
\begin{enumerate}
\item The Goldie closure of $K$ in $M$ is the largest submodule of
$M$ in which $K$ is essential. In particular, $K$ is essential in
its Goldie closure in $M.$
\item The Goldie closure of $K$ in $M$ is the smallest submodule
of $M$ that contains $K$ with no essential extensions in $M$.
\item $K$ is Goldie closed in $M$ if and only if $K$ is a complement
in $M$.
\end{enumerate}
\label{goldie closure}
\end{prop}
This proposition follows from Corollary 7.30 and Proposition 7.44
in \cite{Lam}.

These two propositions have the following  result of R.E. Johnson
(introduced in \cite{John}) as a corollary.

\begin{cor}
Let $R$ be any ring and $M$ a nonsingular $R$-module. There is an
one-to-one correspondence \[\{\mbox{complements in }M\}
\longleftrightarrow\{\mbox{direct summands of }E(M)\}\] given by
$K \mapsto$ the Goldie closure of $K$ in $E(M)$ which is equal to
a copy of $E(K).$ The inverse map is given by $L\mapsto L\cap M.$
\label{Johnson}
\end{cor}

The proof can be found also in \cite{Lam} (Corollary 7.44').

\subsection{Main Results}

Let us consider a finite von Neumann algebra $\A$. Recall that
$\A$ is a nonsingular ring. Since the Goldie torsion theory
coincides with $(\T, \bigP)$ for $\A,$ an $\A$-module $M$ is in
$\bigP$ if and only if it is a nonsingular module.

Johnson's Theorem for the ring $\A$ gives that for every
$\A$-module $M$ in $\bigP$, there is an one-to-one correspondence
\[ \{(\T,\bigP)\mbox{-closed submodules of }M\}
\longleftrightarrow\{\mbox{direct summands of }E(M)\}
\]
given by $K \mapsto$ $\cl_{\T}^{E(M)}(K) = E(K).$ The inverse map
is given by $L\mapsto L\cap M.$ This follows directly from
Johnson's Theorem (Corollary \ref{Johnson}) and Propositions
\ref{Lamova prop o closed submodules} and \ref{goldie closure}.

We shall prove the stronger result in case when the module $M$ is
finitely generated in $\bigP$ (and hence projective). In order to
do that, we need to describe the injective envelope of such $M$.
First we need a lemma.

\begin{lem}
For any $\A$-module $M,$ \[\dim_{\A}(\U\otimes_{\A}M) =
\dim_{\A}(M).\] \label{dim of tenzor product}
\end{lem}
In \cite{Lu5} and \cite{Re} the formula
$\dim_{\U}(\U\otimes_{\A}M) = \dim_{\A}(M)$ is shown. Note the
difference.
\begin{proof} From the short exact sequence
$0\rightarrow\A\rightarrow\U\rightarrow\U/\A\rightarrow0,$ we get
the exact sequence \[0\rightarrow \smallt M \rightarrow M
\rightarrow \U\otimes_{\A}M \rightarrow
\U/\A\otimes_{\A}M\rightarrow0\] since $\smallt M$ is the kernel
of $M\rightarrow\U\otimes_{\A}M$ and $\U$ is $\A$-flat. The module
$\smallt M $ is a submodule of $\T M,$ so it has dimension zero.
To show that the dimensions of $M$ and $\U\otimes_{\A}M$ are the
same, it is sufficient to show that the dimension of
$\U/\A\otimes_{\A}M$ is 0. We shall show a stronger statement: the
module $\U/\A\otimes_{\A}M$ is in $\smallt$ for all $M$, i.e. for
every $a\in \U/\A\otimes_{\A}M $ there is a nonzero-divisor $t$
such that $ta=0.$

Let $a=\sum_{i=1}^n (t_i^{-1}a_i+\A)\otimes_{\A}m_i$ and $t=t_1
t_2\ldots t_n.$ Clearly, $t$ is a nonzero-divisor. Then
$a=\sum_{i=1}^n (t^{-1}t t_i^{-1}a_i+\A)\otimes_{\A}m_i.$ The
fraction $t t_i^{-1}=t_1 t_2\ldots t_n t_i^{-1}$ is in $\A$ by
Lemma \ref{regular_elements_lemma}. Thus, \[ta = \sum_{i=1}^n (t
t_i^{-1}a_i+\A)\otimes_{\A}m_i = \sum_{i=1}^n
(0+\A)\otimes_{\A}m_i=0.\] So, $\U/\A\otimes_{\A}M$ is in
$\smallt.$
\end{proof}

\begin{thm}
\begin{enumerate}
\item If $M$ is an $\A$-module in $\bigP$, then
\[\U\otimes_{\A}M\subseteq E(M).\]
\item If $M$ is a finitely generated projective $\A$-module, then
\[\U\otimes_{\A}M = E(M).\]
\end{enumerate}
\label{envelope is tensor}
\end{thm}
\begin{proof}
(1) Let $M$ be in $\bigP.$ Then, $M$ is flat. Hence $0=\smallt M =
\tor_1^{\A}(\U/\A,M)$, so $M$ embeds in $\U\otimes_{\A}M.$ By
previous lemma, $M$ and $\U\otimes_{\A}M$ have the same dimension,
so the closure of $M$ in $\U\otimes_{\A}M$ with respect to $(\T,
\bigP)$ is equal to entire $\U\otimes_{\A}M.$ Since a submodule of
a nonsingular module is an essential submodule of its closure
(Proposition \ref{goldie closure}), we have
$M\subseteq_e\U\otimes_{\A}M.$ The injective envelope $E(M)$ of
$M$ is the maximal essential extension of $M$ and so
$\U\otimes_{\A}M$ is contained in a copy of $E(M).$

(2) From (1), we have that $M\subseteq_e\U\otimes_{\A}M.$ So, the
injective envelopes of $M$ and $\U\otimes_{\A}M$ are the same. To
show the claim, it is sufficient to show that $\U\otimes_{\A}M$ is
an injective $\A$-module.

Since $M$ is a finitely generated projective module, there is a
positive integer $n$ and a module $N$ such that $M\oplus N =
\A^n.$ So, $\U\otimes_{\A}M$ is a direct summand of
$\U\otimes_{\A}\A^n \cong \U^n.$ Since $\U$ is $\A$-injective,
$\U^n$ is $\A$-injective as well and, so is its direct summand
$\U\otimes_{\A}M.$
\end{proof}

The following is an example of a nonsingular $\A$-module with
strict inclusion in part (1) of the above theorem.

\begin{exmp} Consider the group von Neumann algebra $\NZ$ of the group
$\Zset.$ Example 8.34 in \cite{Lu5} produces an example of a
$\NZ$-ideal $I$ such that $M=\NZ/I$ is a flat module of dimension
zero (so this proves that $\smallt\subsetneq\T$). The short exact
sequence $0\rightarrow I\rightarrow \NZ\rightarrow M\rightarrow0$
gives us
\[0\rightarrow \UZ\otimes_{\NZ}I\rightarrow\UZ\otimes_{\NZ}\NZ = \UZ
\rightarrow\UZ\otimes_{\NZ}M\rightarrow0.\] Since $I$ and $M$ are
flat modules, the modules $\UZ\otimes_{\NZ}I$ and
$\UZ\otimes_{\NZ}M$ are nonzero (they would vanish just if $I$ and
$M$ were in $\smallt$). So, $0 \neq \UZ\otimes_{\NZ} I \subsetneq
\UZ\otimes_{\NZ}\NZ = \UZ.$ Since the dimension of $M$ is zero,
$I\subseteq_e\NZ.$ Thus, $E(I) = E(\NZ) =\UZ.$ Hence,
$\UZ\otimes_{\NZ}I\subsetneq E(I).$
\end{exmp}

As a corollary of Johnson's result (Corollary \ref{Johnson}) and
previous theorem, we obtain the following.

\begin{cor}
For every finitely generated projective $\A$-module $M$, there is
an one-to-one correspondence \[\{\mbox{direct summands of }M\}
\longleftrightarrow\{\mbox{direct summands of }E(M)\}\] given by
$K \mapsto$ $\U\otimes_{\A}K = E(K).$ The inverse map is given by
$L\mapsto L\cap M.$ \label{direct summands are isomorphic}
\end{cor}
\begin{proof}
Recall that the Johnson's result gives us the correspondence
\[ \{(\T,\bigP)\mbox{-closed submodules
of }M\} \longleftrightarrow\{\mbox{direct summands of }E(M)\} \]
with $K \mapsto E(K)$ and the inverse $L\mapsto L\cap M.$ Since
$M$ is finitely generated projective, the previous theorem gives
us that $\U\otimes_{\A}K = E(K)$ if $K$ is a direct summand of
$M.$ Thus, it is sufficient to prove that a submodule $K$ of $M$
is a direct summand of $M$ if and only if it is
$(\T,\bigP)$-closed. Clearly, if $K$ is a direct summand of $M$
then $M/K$ is projective so $0 = \T(M/K) = \cl_{\T}(K)/K.$
Conversely, if $0 = \cl_{\T}(K)/K = \T(M/K),$ then $M/K$ is a
finitely generated module in $\bigP$ and, hence, projective. Thus,
$K$ is a direct summand of $M.$
\end{proof}

Using this result, we can obtain the explicit description of the
map $\mu^{-1}:$ Proj($\U$) $\rightarrow$ Proj($\A$) that induces
the inverse of the isomorphism $\mu: K_0(\A)\rightarrow K_0(\U).$

\begin{thm} There is an one-to-one correspondence between
Goldie closed ideals of $\A$ and direct summands of $\U$ given by
$I \mapsto$ $\U\otimes_{\A}I = E(I).$ The inverse map is given by
$L\mapsto L\cap \A.$ This correspondence induces an isomorphism of
monoids $\mu :\mathrm{Proj}(\A)$ $\rightarrow$ $\mathrm{Proj}(\U)$
and an isomorphism \[\mu: K_0(\A)\rightarrow K_0(\U)\] given by
$[P]\mapsto[\U\otimes_{\A}P]$ with the inverse
$[Q]\mapsto[Q\cap\A^n]$ if $Q$ is a direct summand of $\U^n.$
\label{moja K 0 theorem}
\end{thm}
The proof follows directly from Corollary \ref{direct summands are
isomorphic}.

\section{$\smallt-\T\p-\bigP$ FILTRATION AND THE CAPACITY}
\label{capacity}

In this section, we describe the three parts $\smallt M$, $\T\p M$
and $\bigP M $ of an $\A$-module $M$ as the certain submodules of
a free cover of $M$. Then, we prove a formula that gives the
capacity of an $\A$-module via the capacity of these submodules.

Let us begin with a technical lemma that we need.

\begin{lem}
Let $F$ be a finitely generated free (or projective) $\A$-module
and $K$ its submodule. Let $K_i,$ $i\in I,$ be any directed family
of finitely generated submodules of $K$ (directed with respect to
the inclusion maps) such that the directed union $\dirlim K_i$ is
equal to $K.$ Then
\[\cl_{\smallt}(K) = \dirlim \cl_{\smallt}(K_i) =
\dirlim\cl_{\T}(K_i).\] \label{cl t K VS dirlim}
\end{lem}
\begin{proof}
First note that $\cl_{\T}(K_i) = \cl_{\smallt}(K_i)$ because
modules $K_i$ are projective as finitely generated submodules of
the projective module $F$ ($\A$ is semihereditary). So $F/K_i$ are
finitely presented modules. Since $\T = \smallt$ for the class of
finitely presented modules, we have that $\cl_{\T}(K_i)/K_i =
\T(F/K_i) = \smallt(F/K_i) = \cl_{\smallt}(K_i)/K_i,$ so the
second equality follows. Let us show now that $\cl_{\smallt}(K) =
\dirlim\cl_{\T}(K_i).$

Since $K_i\subseteq K,$ we have that $\cl_{\T}(K_i) =
\cl_{\smallt}(K_i)\subseteq \cl_{\smallt}(K).$ So $\dirlim
\cl_{\T}(K_i)\subseteq \cl_{\smallt}(K).$

For the converse, look at the quotient $Q =
\cl_{\smallt}(K)/\dirlim \cl_{\T}(K_i).$ We shall show it is equal
to 0. By applying direct limit functor to the short exact
sequence:
\[0\rightarrow\cl_{\T}(K_i)\rightarrow\cl_{\smallt}(K)\rightarrow
\cl_{\smallt}(K)/\cl_{\T}(K_i) \rightarrow 0,\] we have that the
direct limit $P$ of the quotients $\cl_{\smallt}(K)/\cl_{\T}(K_i)$
is isomorphic to $Q$. We shall show that $Q$ is trivial by showing
that $P$ is flat and $Q$ is in $\smallt.$ This will give us that
$P\cong Q$ is both in $\smallt$ and $\p.$ So, it must be zero.

To show that $P$ is flat, note that the module
$\cl_{\smallt}(K)/\cl_{\T}(K_i)$ is a submodule of the module
$F/\cl_{\T}(K_i)$ for every $i$ in $I.$ The latter module is
finitely generated and projective, so it is flat. Since a
submodule of a flat $\A$-module is flat ($\A$ is semihereditary),
$\cl_{\smallt}(K)/\cl_{\T}(K_i)$ is flat. Since direct limits
preserves flatness, the module $P =
\dirlim\cl_{\smallt}(K)/\cl_{\T}(K_i)$ is flat.

On the other hand, $K = \dirlim K_i \subseteq \dirlim
\cl_{\T}(K_i)$ so the module $Q$ is a quotient of the module $
\cl_{\smallt}(K)/K = \smallt(F/K)$ which is in $\smallt.$ Thus,
$Q$ is in $\smallt.$ This finishes the proof.
\end{proof}

The next proposition describes the three parts $\smallt M$, $\T\p
M$ and $\bigP M $ of an $\A$-module $M$ via certain submodules of
a free cover of $M$.

\begin{prop}
Let $M$ be a finitely generated $\A$-module. Let $F$ be a finitely
generated free (or projective) module that maps onto $M$ by some
map $f$. Let $K$ be the kernel of $f$. Let $K_i,$ $i\in I,$ be any
directed family of finitely generated submodules of $K$ (directed
with respect to the inclusion maps) such that the union $\dirlim
K_i$ is equal to $K.$ Then
\begin{enumerate}
\item \[\smallt M =\cl_{\smallt}(K)/K \cong
\lim_{\overrightarrow{i\in I}}\left(\cl_{\T}(K_i)/K_i\right)=
\lim_{\overrightarrow{i\in I}}\T(F/K_i).\] $M$ is flat if and only
if $\cl_{\smallt}(K)=\underrightarrow{\lim}\; \cl_{\T}(K_i) = K.$
\item \[\T\p M
\cong \cl_{\T}(K)/\cl_{\smallt}(K)
 \cong \lim_{\overrightarrow{i\in
 I}}\left(\cl_{\T}(K)/\cl_{\T}(K_i)\right).\]

\item \[\bigP M \cong F/\cl_{\T}(K).\]

All closures are taken in $F.$
\end{enumerate}
\label{t-Tp-P}
\end{prop}
\begin{proof}
The first and the last equality in (1) follow since $\smallt M =
\smallt(F/K) =\cl_{\smallt}(K)/K$ and $\T(F/K_i) =
\cl_{\T}(K_i)/K_i$ by Proposition \ref{properties of closure}. We
have the middle isomorphism in (1) because $\cl_{\smallt}(K) =
\dirlim \cl_{\T}(K_i)$ (Lemma \ref{cl t K VS dirlim}) and because
the direct limit functor is exact.

$M$ is flat iff $\smallt M=0$ iff $\cl_{\smallt}(K)=K.$ This is
equivalent with $\cl_{\smallt}(K) = \dirlim \cl_{\T}(K_i)$ by
Lemma \ref{cl t K VS dirlim}.

From the proof of part (4) of Proposition \ref{odnosi Tova i
Pova}, we have that $\T\p M = \T M/\smallt M.$ $ \T M/\smallt M =$
$\T(F/K)/\smallt(F/K) =$
$(\cl_{\T}(K)/K)/(\cl_{\smallt}(K)/K)\cong$ $
\cl_{\T}(K)/\cl_{\smallt}(K).$ The second isomorphism in (2)
follows by Lemma \ref{cl t K VS dirlim} and by exactness of the
direct limit functor.

Part (3) follows from Proposition \ref{properties of closure}:
$\bigP M = \bigP (F/K) \cong F/\cl_{\T}(K).$
\end{proof}

\subsection{Capacity}

We have seen that every $\A$-module $M$ consists of three parts: a
cofinal-measurable $\smallt$-part, a flat, zero-dimensional
$\T\p$-part and a $\bigP$-part. The dimension measures the
$\bigP$-part faithfully. The cofinal-measurable part can also be
measured faithfully. The invariant that measures it is the {\em
Novikov-Shubin invariant}, $\alpha(M).$ Sometimes, for
convenience, the reciprocal of the $\alpha(M)$ is considered. The
reciprocal $c(M) = \frac{1}{\alpha(M)}$ is called the {\em
capacity}.

The Novikov-Shubin invariant is defined first for a finitely
presented $\A$-module $M$ and then the definition is extended to
every $\A$-module. First, the two finitely generated projective
modules $P_0$ and $P_1$ with quotient $M$ are considered:
$0\rightarrow P_1\rightarrow P_0\rightarrow M\rightarrow 0$. Then,
the equivalence $\nu$ of the category of finitely generated
projective $\A$-modules and the finitely generated Hilbert
$\A$-modules from \cite{Lu1} is used to get the morphism $\nu(i):
\nu(P_1)\hookrightarrow\nu(P_0)$ of finitely generated Hilbert
$\A$-modules. We define the Novikov-Shubin invariant $\alpha(M)$
of $M$ via the Novikov-Shubin invariant of the morphism $\nu(i).$

Let $f: U\rightarrow V$ be a morphism of two finitely generated
Hilbert $G$-modules. Then the operator $f^*f$ is positive. Let
$\{\;E_{\lambda}^{f^*f}\;|\:\lambda\in \Rset\;\}$ be the family of
spectral projections of $f^*f.$ Define the {\em spectral density
function} of $f$ by \[F(f):[0,
\infty)\rightarrow[0,\infty],\mathrm{    }\;\;\;
\lambda\mapsto\dim_{\A}(\im
(E_{\lambda^2}^{f^*f}))=\tr(E_{\lambda^2}^{f^*f}).\]

The {\em Novikov-Shubin invariant of a morphism} $f:U\rightarrow
V$ of finitely generated Hilbert $G$-modules by
\[\alpha (f) = \liminf_{\lambda\rightarrow
0^+}\frac{\ln(F(f)(\lambda) - F(f)(0))}{\ln \lambda}\in[0,
\infty],\] if $F(f)(\lambda) > F(f)(0)$ for all $\lambda>0.$ If
not, we let $\alpha (f) = \infty^+$ where $\infty^+$ is a new
symbol. We define an ordering on the set
$[0,\infty]\cup\{\infty^+\}$ by the standard ordering on $\Rset$
and $x<\infty<\infty^+$ for all $x\in\Rset.$

If $M$ is a finitely presented module with finitely generated
projective modules $P_0$ and $P_1$ and the short exact sequence
$0\rightarrow P_1\rightarrow P_0\rightarrow M\rightarrow 0,$ the
Novikov-Shubin invariant $\alpha(M)$ measures the $\smallt =
\T$-part of $M$: smaller $\alpha(M)$ corresponds to a larger
difference between $P_0$ and its closure in $P_1$, i.e. to larger
$\T M=\smallt M$. $\alpha(M)$ is $\infty^+$ if and only if $M$ is
projective itself, i.e. $\cl_{\T}^{P_1}(P_0)/P_0=\T M=\smallt M=0$
(for details see \cite{Wegner}).

The capacity $c(M)\in \{0^-\}\cup[0,\infty]$ of such finitely
presented $M$ is $c(M) = \frac{1}{\alpha (M)},$ where
$0^-:=(\infty^+)^{-1},$ $0^{-1}=\infty$ and $\infty^{-1}=0.$

Next, we define the capacity of measurable module $M$ (quotients
of finitely presented $\T$-modules) as follows
\[
c(M) =\inf\;\{\; c(L)\;|\;L\mbox{ fin. presented,
zero-dimensional, }M\mbox{ quotient of }L\;\}.\]

Finally, the capacity of arbitrary $\A$-module $M$ is defined as
\[ c(M) =\sup\;\{\; c(N)\;|\;N\mbox{ measurable submodule of
}M\;\}.\]

The following proposition shows that the capacity measures
faithfully $\smallt$-part of any $\A$-module. Also, we can use
capacity to check if an $\A$-module is flat.

\begin{prop}
Let $M$ be an $\A$-module. Then
\begin{center}
\begin{enumerate}
\item $c(M) = c(\smallt M)$ and $c(\p M) = 0^-.$
\item $c(\smallt M) = 0^-$ if and only if $\smallt M=0.$
\item $M$ is flat if and only if $c(M) = 0^-.$
\end{enumerate}
\end{center}
\label{prop about capacity}
\end{prop}
\begin{proof} (1) Since any measurable submodule of $M$ is in
$\smallt M,$ we have that $c(M)\leq c(\smallt M).$ The converse
clearly holds so $c(M) = c(\smallt M).$ Since $\smallt\p M = 0,$ $
c(\p M) = c(\smallt\p M) = c(0) = 0^-.$

(2) Clearly $c(0)=0^-$. Since $c(M)=c(\smallt M)$ by (1), it is
sufficient to show (2) for a module $M$ in $\smallt$. If $M$ is
finitely presented and in $\smallt$, then $M = \smallt M= \T M =
0$ iff $c(M)=0^-$ by the remarks following the definition of
Novikov-Shubin invariant.

If $M$ is a measurable module with capacity $0^-$, there is a
finitely presented, zero-dimensional module $L$ with capacity
$0^-$ which has $M$ as a quotient. But then such $L$ must be zero
by the previous case, and so $M=0$ as well.

Now, let $M$ be any module in $\smallt$ with capacity $0^-$. Then
every measurable submodule of $M$ has capacity $0^-.$ But, by
previous case, that means that every measurable submodule of $M$
is 0. Then $M$ has to be 0 as well because a cofinal-measurable
module is the directed union of its measurable submodules.

(3) $M$ is flat  $\Leftrightarrow M=\p M$ $\Leftrightarrow \smallt
M = 0$ $\Leftrightarrow c(\smallt M)= 0^-$ $\Leftrightarrow
c(M)=0^-.$
\end{proof}

The capacity also has the following properties.
\begin{prop}
\begin{enumerate}
\item If $0\rightarrow M_0\rightarrow M_1\rightarrow M_2\rightarrow
0$ is a short exact sequence of $\A$-modules, then
\begin{itemize}
\item[i)] $c(M_0) \leq c(M_1);$
\item[ii)] $c(M_2) \leq c(M_1)$ if $M_1$ is in $\smallt;$
\item[iii)] $c(M_1) \leq c(M_0) + c(M_2)$ if $M_1$ is in $\T.$
\end{itemize}

\item If $M = \bigcup_{i \in I}M_i$ is a directed union, then
$c(M) = \sup\{\;c(M_i)\; |\; i\in I\;\}.$

\item If $M = \bigoplus_{i\in I} M_i,$ then $c(M) =
\sup\{\;c(M_i)\;|\;i \in I\;\}.$
\item If $M = \dirlim_{i\in I} M_i$ is a direct limit of a directed
system with structure maps $f_{ij}: M_i\rightarrow M_j,$ $i\leq
j,$ then \[ c(M) \leq \lim\inf\{\;c(M_i)\;|\;i\in I\}.\] If $M_i$
is measurable for all $i\in I$ and the maps $M_i\rightarrow M_j$
are surjective for all $i, j\in I$ such that $i\leq j,$ then $c(M)
= \inf\{\;c(M_i)\;|\;i\in I\;\}. $
\end{enumerate}
\label{prop of c}
\end{prop}

For the proof, see \cite{Lu4} or \cite{Wegner}.

In \cite{Lu4}, the following formula for computing the capacity of
a measurable module is given.

\begin{lem}
Let $M$ be a measurable $\A$-module. Let $f: F\rightarrow M$ be a
surjection of a finitely generated free (or projective) module $F$
onto $M.$ Then
\[ c(M) =\inf\;\{\; c(F/K)\;|\;K\subseteq\ker
f\mbox{ fin. gen. and } \dim_{\A}(K) = \dim_{\A}(F)\;\}.\] The set
on the right hand side is nonempty if and only if $M$ is
measurable. \label{formula_for_c_Luck}
\end{lem}

We shall prove a more general formula for capacity. The formula
will show that we can use the modules $K_i,$ $i\in I$ from the
setting like the one in Proposition \ref{t-Tp-P} to calculate the
capacity of a finitely generated module.

\begin{prop}
Let $M$ be a finitely generated $\A$-module. Let $F$ be a finitely
generated free (or projective) module that maps onto $M$ by some
map $f$. Let $K$ be the kernel of $f$. Let $K_i,$ $i\in I,$ be any
directed family of finitely generated submodules of $K$ (directed
with respect to the inclusion maps) such that the union $\dirlim
K_i$ is equal to $K.$ Then
\begin{eqnarray}
c(M) & = & \sup_{i\in I}\;
c\left(\cl_{\T}(K_i)/\cl_{\T}^{K}(K_i)\right)
  =  \sup_{i\in I}\inf_{j\geq i}\;
  c\left(\cl_{\T}(K_i)/\cl_{\T}^{K_j}(K_i)\right)\nonumber\\
& = & \sup_{i\in I}\inf_{j\geq i}\; c\left(\cl_{\T}(K_i)/K_j\cap
\cl_{\T}(K_i)\right). \nonumber
\end{eqnarray} \label{my capacity}
\end{prop}
\begin{proof}
Let us first note that $K_j\cap \cl_{\T}(K_i) =
\cl_{\T}^{K_j}(K_i)$ by part (6) of Proposition \ref{properties of
closure}. Thus, the third equality follows.

Recall that $c(M) = c(\smallt M)$ and $\smallt M =
\dirlim\cl_{\T}(K_i)/K_i$ (by Proposition \ref{t-Tp-P}). So,
$\smallt M$ is the directed union of the images of maps $f_i:
\cl_{\T}(K_i)/K_i\rightarrow\cl_{\smallt}(K)/K.$ The kernel of
$f_i$ is $\cl_{\smallt}^K(K_i)/K_i$ and so the image of $f_i$ is
isomorphic to the quotient $\cl_{\T}(K_i)/\cl_{\smallt}^K(K_i)
=\cl_{\T}(K_i)/K\cap \cl_{\smallt}(K_i)= \cl_{\T}(K_i)/K\cap
\cl_{\T}(K_i).$ Thus, $\smallt M \cong \bigcup_{i\in
I}\;\cl_{\T}(K_i)/K\cap \cl_{\T}(K_i)$ and so the first equality
follows by part (2) of Proposition \ref{prop of c}.

Now, let us fix $i$ in $I$ and look at the quotient
$\cl_{\T}(K_i)/K\cap \cl_{\T}(K_i)$ again. Since $K$ is the
directed union of $K_j$ where $j\geq i,$ we have that
\[\cl_{\T}(K_i)/K\cap \cl_{\T}(K_i) = \dirlim_{j\geq i}
\cl_{\T}(K_i)/K_j\cap \cl_{\T}(K_i)\] by the exactness of the
direct limit functor. If $k\geq j\geq i$ the structure map
$\cl_{\T}(K_i)/K_j\cap
\cl_{\T}(K_i)\rightarrow\cl_{\T}(K_i)/K_k\cap \cl_{\T}(K_i)$ is
onto. Moreover, the module $\cl_{\T}(K_i)/K_j\cap \cl_{\T}(K_i)$
is measurable for every $j\geq i$ since it is a quotient of the
finitely presented zero-dimensional module $\cl_{\T}(K_i)/K_i.$

Thus, the two conditions of the second part of (4) in Proposition
\ref{prop of c} are satisfied and, we obtain that the capacity of
the quotient $\cl_{\T}(K_i)/K\cap \cl_{\T}(K_i)$ is equal to the
infimum of the capacities of $\cl_{\T}(K_i)/K_j\cap \cl_{\T}(K_i)$
for $j\geq i.$ This gives us the second equality.
\end{proof}

The formula from the above theorem agrees with the condition from
(1) in Proposition \ref{t-Tp-P}: $M$ is flat iff $0^- = c(M)$ iff
$\dirlim \cl_{\T}(K_i) = K.$ Taking $K_i'=\cl_{\T}(K_i)$, we
obtain a family $\{K_i'\}_{i\in I}$ such that
$K_j'\cap\cl_{\T}(K_i') = K_j'\cap K_i' = K_i' = \cl_{\T}(K_i')$
for all $j\geq i.$ So, the quotient
$\cl_{\T}(K_i)/K_j\cap\cl_{\T}(K_i)$ is equal to zero.

\section{INDUCTION}
\label{Induction}

Let $H$ be a Hilbert space and $\A$ a von Neumann algebra in
$\B(H).$ A $C^*$-subalgebra $\B$ of $\A$ is a {\em von Neumann
subalgebra} of $\A$ if $B''=B$ where the commutants are computed
in $\B(H)$ (equivalently $\B$ is closed with respect to weak or
strong operator topology). If $\A$ is finite with normal and
faithful trace $\tr_{\A},$ the restriction $\tr_{\B}=\tr_{\A}|\B$
of $\tr_{\A}$ to $\B$ is a normal and faithful trace on $\B,$ so
$\B$ is finite as well. If $\B$ is a von Neumann subalgebra of a
finite von Neumann algebra $\A$, the only normal and faithful
trace on $\B$ that we consider is the restriction of the normal
and faithful trace on $\A.$ Note that algebra $\B$ might have
other normal and faithful trace functions besides $\tr_{\A}|\B.$

If $\B$ is a von Neumann subalgebra of a finite von Neumann
algebra $\A,$ and $M$ is a $\B$-module, we define the {\em
induction of $M$} as the $\A$-module
\[i_*(M) = \A\otimes_{\B}M.\]

In this section, we shall prove that the dimension and capacity
are both preserved by induction.

In case when $\A=\vng$ is a group von Neumann algebra and
$\B=\vnh$ where $H$ is a subgroup of $G,$ the result that the
dimension is preserved under induction is given in \cite{Lu2}. We
shall show that the same holds for finite von Neumann algebras.
The inequality $c(M)\leq c(i_*(M))$ for a group von Neumann
algebra case is proven in \cite{Wegner}. We shall prove that the
equality
\[c(M)= c(i_*(M))\] holds for any finite von Neumann algebra.

Let $\B$ be a von Neumann subalgebra of a finite von Neumann
algebra $\A.$ Clearly, $i_*(\B) = \A.$ From the definition of
$i_*$, it follows that $i_*$ is a covariant functor from the
category of $\B$-modules to the category of $\A$-modules which
maps a direct sum to a direct sum, a finitely generated module to
a finitely generated module and a projective module to a
projective module. Also, $i_*$ commutes with direct limits.

First we shall prove the generalization of Theorem 3.3 from
\cite{Lu2}.
\begin{prop}
Let $\B$ be a von Neumann subalgebra of a finite von Neumann
algebra $\A.$ The induction $i_*$ is a faithfully flat functor
from the category of $\B$-modules to the category of $\A$-modules.
If $M$ is a $\B$-module, then \[\dim_{\B}(M) =
\dim_{\A}(i_*(M)).\]
\end{prop}
\begin{proof}
Theorem 3.3 from \cite{Lu2} states the same about the functor
$i_*$ but just in the case of a group von Neumann algebra $\vng$
and its subalgebra $\vnh$ where $H$ is a subgroup of $G.$ L\"uck's
proof of Theorem 3.3 in \cite{Lu2} consists of seven steps. In
step 1, it is shown that $i_*$ preserves the dimension of a
finitely generated projective module. In step 2, it is shown that
this is true for finitely presented modules as well and that
$\tor^{\vnh}_1(\vng, M)=0$ if $M$ is a finitely presented
$\vnh$-module. Step 3 shows that $\tor^{\vnh}_1(\vng, M)=0$ if $M$
is finitely generated. In step 4, L\"uck shows that $i_*$ is an
exact functor. Steps 5 and 6 prove that $i_*$ preserves the
dimension. Finally, in step 7 L\"uck shows that $i_*$ is faithful.

To prove this more general theorem about finite von Neumann
algebras, the only modification to the proof of L\" uck's Theorem
3.3 must be made in the first two steps -- the remaining steps of
the proof hold for any finite von Neumann algebra without any
modifications.

In L\"uck's proof of step 1, the key observation was that the
standard trace on $\vnh$ is the restriction of the standard trace
on $\vng$ so $\tr_{\vnh}(a)=\tr_{\vng}(i(a))$ where $i$ is the
inclusion $\vnh\rightarrow\vng.$ This observation remains true for
finite von Neumann algebras as well since $\tr_{\B}=\tr_{\A}|\B.$
With this in mind, L\"uck's proof of step 1 holds for finite von
Neumann algebras.

Before proving step 2, let us note that $\A\otimes_{\B}\lb$ is a
dense subspace of $\la.$ This is the case we can identify $\A$
with $\A\otimes_{\B}\B$ and $\la$ with $l^2(\A\otimes_{\B}\B).$
$\A\otimes_{\B}\lb$ is dense in $\la\otimes_{\B}\lb$ which is
dense in $l^2(\A\otimes_{\B}\B).$

Now, let $M$ be a finitely presented module. $M$ is the direct sum
of finitely generated projective module $\bigP M$ and
zero-dimensional $\T M$. By step 1, step 2 clearly holds for
$\bigP M.$ So, it is sufficient to consider finitely presented
modules in $\T.$ If $M$ is such module, there is a nonnegative
integer n and an injective map $f:\B^n\rightarrow\B^n$ with
$f^*=f$ and $M=$ Coker $f$ (see \cite{Lu1} for the proof of this
fact). In order to prove step 2, it is sufficient to prove that
$i_*(f): \A^n\rightarrow\A^n$ is injective with the cokernel of
dimension zero.

The functor $\nu$ from Theorem \ref{equivalence} is weak exact
(see \cite{Lu1} for proof). Thus, the image of $\nu(f):
\lb^n\rightarrow\lb^n$ is dense in $\lb^n.$ Since the diagram
\[
\begin{array}{cccc}
i_*(\nu(f)): & \A\otimes_{\B}\lb^n & \rightarrow &
\A\otimes_{\B}\lb^n \\ &\downarrow & &\downarrow \\
\nu(i_*(f)): & l^2(\A\otimes_{\B}\B^n) &
\rightarrow & l^2(\A\otimes_{\B}\B^n) \\
\end{array}
\]
commutes and $\A\otimes_{\B}\lb^n$ is dense in
$l^2(\A\otimes_{\B}\B^n),$ the image of the map $\nu(i_*(f))$ is
dense. The functor $\nu$ is such that $\nu(f^*)=(\nu(f))^*$ (see
\cite{Lu1}). Thus, $\nu(i_*(f))$ is selfadjoint since $f$ is. The
image of $\nu(i_*(f))$ is dense and so the kernel of $\nu(i_*(f))$
is trivial. Since $\nu^{-1}$ is exact (see \cite{Lu1}), the kernel
of $i_*(f)$ is trivial also. So,
$0\rightarrow\A^n\rightarrow\A^n\rightarrow \mathrm{Coker}(
i_*(f))\rightarrow0.$ Hence, $\dim_{\A}(\mathrm{Coker} (i_*(f)))=
0$ by the additivity of the dimension function (Proposition
\ref{propofdim}). This finishes the proof of step 2 for the case
of finite von Neumann algebras.
\end{proof}

We can define the induction functor on the category of finitely
generated Hilbert $B$-modules as follows. Let $U$ be a finitely
generated Hilbert $\B$-module endowed with a $\B$-valued inner
product $\langle\;\;,\;\;\rangle_{\B}$. Then, $\A\otimes_{\B} U$
has a pre-Hilbert structure via the $\A$-valued inner product:
\[\langle a_1\otimes u_1 , a_2\otimes u_2\rangle = a_1\;\langle u_1,
u_2\rangle_{\B}\;a_2^*\] composed with $\tr_{\A}.$

Define the induction $i_*(U)$ to be the completion of the
pre-Hilbert space $\A\otimes_{\B} U$. The induction of a morphism
$g:U\rightarrow V $ of finitely generated Hilbert $B$-modules is
the induced map $i_*(U)\rightarrow i_*(V).$ With this definition,
the above commutative diagram gives us that $\nu\circ i_* =
i_*\circ \nu$ on the category of finitely generated projective
$\B$-modules.

Now we shall prove the generalization of Proposition 4.4.10 from
\cite{Wegner}.
\begin{prop}
Let $\B$ be a von Neumann subalgebra of a finite von Neumann
algebra $\A.$ Then,
\begin{enumerate}
\item If $M$ is a $\B$-module in $\smallt$, then $i_*(M)$ is also
in $\smallt$ and $c(M) = c(i_*(M)).$
\item If $M$ is any $\B$-module, $c(M)\leq c(i_*(M)).$
\end{enumerate}
\label{i tudja}
\end{prop}
\begin{proof}
In \cite{Wegner}, Wegner proves the proposition about the group
von Neumann algebras in fours steps. In the first step, he proves
that (1) holds for finitely presented modules in $\T$. In the
second step, he proves that (1) holds for measurable modules and
in the third that (1) holds for all modules in $\smallt$. In the
fourth step, he proves (2).

To prove this more general theorem about finite von Neumann
algebras, the only modification to the Wegner's proof must be made
in the first step -- the remaining steps of the proof hold for any
finite von Neumann algebra without any modifications.

Let $M$ be a finitely presented $\B$-module in $\smallt$. Then
there is a short exact sequence $0\rightarrow\B^n\rightarrow\B^n
\rightarrow M\rightarrow0$ for some nonnegative integer $n$. Since
$i_*$ is exact, $i_*(M)$ is also finitely presented and in
$\smallt.$ The capacity of $M$ is defined as the capacity of the
map $\nu(f): \lb^n\rightarrow\lb^n$ and the capacity of $i_*(M)$
as the capacity of the map $\nu(i_*(f)): \la^n\rightarrow\la^n.$
To show that $c(\nu(f))=c(\nu(i_*(f))),$ it is sufficient to show
that $c(g) = c(i_*(g))$ for every morphism of finitely generated
Hilbert $\B$-modules. If $g$ is such a morphism, it is easy to
check that the spectral projections satisfy that
$i_*(E_{\lambda}^{g^*g}) = E_{\lambda}^{i_*(g^*g)}.$ This gives us
that g and $i_*(g)$ have the same spectral density functions
$F(g)= F(i_*(g))$ and so $c(g) = c(i_*(g)).$

The remainder of Wegner's proof holds for finite von Neumann
algebras without any modifications.
\end{proof}

Now, we shall prove that the equality $c(M) = c(i_*(M))$ holds for
{\em all} $\B$-modules $M$. To prove that, we shall use part (1)
of the following proposition.

\begin{prop}
Let $\B$ be a von Neumann subalgebra of a finite von Neumann
algebra $\A.$ Then
\begin{enumerate}
\item $i_*(\smallt M)  =  \smallt i_*(M)$ and $i_*(\p M)  =  \p
i_*(M).$
\item $i_*(\T M) = \T i_*(M)$  and $i_*(\bigP M)  = \bigP i_*(M).$
\end{enumerate}
\label{i commutes}
\end{prop}
\begin{proof} (1) If $M$ is flat (i.e. in $\p$), then it is a direct
limit of finitely generated projective modules. But since $i_*$
preserves both the direct limits and finitely projective modules,
the module $i_*(M)$ is also a direct limit of finitely generated
projective modules and, hence, flat.

By \ref{i tudja}, $i_*(\smallt M)$ is in $\smallt$. Since, $i_*$
is exact, we have the short exact sequence
\[0\rightarrow i_*(\smallt M)\rightarrow i_*( M)
\rightarrow i_*(\p M)\rightarrow0\] where $i_*(\smallt M)$ is in
$\smallt$ and $i_*(\p M)$ is in $\p.$ But that means that
$i_*(\smallt M)  =  \smallt i_*(M)$ and $i_*(\p M) = \p i_*(M).$

(2) Let $M$ be in $\bigP.$ Recall that a module is in $\bigP$ if
and only if every finitely generated submodule is projective. So,
$M$ is equal to the directed union of its finitely generated and
projective submodules $M_i,$ $i\in I.$ Then $i_*(M)$ is the
directed union of finitely generated projective modules
$i_*(M_i).$ Let $N$ be a finitely generated submodule of $i_*(M).$
Then there is $i\in I$ such that $N$ is contained in $i_*(M_i).$
But, since $\A$ is semihereditary, we have that $N$ is projective
as well. So, $i_*(M)$ is in $\bigP.$

If $M$ is in $\T,$ then $i_*(M)$ is also in $\T$ since
$\dim_{\B}(M) = \dim_{\A}(i_*(M)).$ Since $i_*$ preserves both
$\T$ and $\bigP$ and $i_*$ is exact, we have that $i_*(\T M)=\T
i_*(M)$ and $i_*(\bigP M) =\bigP i_*(M).$
\end{proof}

\begin{thm} Let $\B$ be a von Neumann subalgebra of a finite von
Neumann algebra $\A.$ Then \[c(M)= c(i_*(M)).\] \label{improved
formula for c}
\end{thm}
\begin{proof}
\[
\begin{array}{rcll}
c(M) & = & c(\smallt M) &
\mbox{(by Proposition \ref{prop about capacity})}\\
 & = & c(i_*(\smallt M)) & \mbox{(by Proposition \ref{i tudja})}\\
 & = & c(\smallt i_*(M)) &
 \mbox{(by Proposition\ref{i commutes})}\\
 & = & c(i_*(M)) &
 \mbox{(by Proposition \ref{prop about capacity})}. \\
\end{array}
\]
\end{proof}

\end{document}